\documentclass{article}
\usepackage{amssymb,amsmath}
\usepackage{graphics}

\def\qed{\hfill {\hbox{${\vcenter{\vbox{               
   \hrule height 0.4pt\hbox{\vrule width 0.4pt height 6pt
   \kern5pt\vrule width 0.4pt}\hrule height 0.4pt}}}$}}}
  
\def\tr{\triangleright}

\newtheorem{theorem}{Theorem}
\newtheorem{definition}{Definition}
\newtheorem{lemma}[theorem]{Lemma}
\newtheorem{proposition}[theorem]{Proposition}
\newtheorem{corollary}[theorem]{Corollary}
\newtheorem{example}{Example}
\newtheorem{remark}{Remark}

\newenvironment{proof}[1][Proof]{\medskip\noindent{\bf #1.}\quad}%
{\qed\par\medskip}

\author{{\begin{tabular}{c} Sam Nelson \\
\small{\texttt{knots@esotericka.org}}\end{tabular}}
\and {\begin{tabular}{c}Chau-Yim Wong \\
\small{\texttt{jason@math.ucr.edu}}\end{tabular}}
\and
\small{\begin{tabular}{c} Department of Mathematics, 
University of California, Riverside, \\ 900 University Avenue,
Riverside, CA, 92521\end{tabular} }}

\date{}

\title{\Large \textbf{On the Orbit Decomposition of Finite Quandles}}

\begin{document}

\maketitle

\begin{abstract}
We study the structure of finite quandles in terms of subquandles. Every 
finite quandle $Q$ decomposes in a natural way as a union of disjoint 
$Q$-complemented subquandles; this decomposition coincides with the usual
orbit decomposition of $Q$. Conversely, the structure of a finite quandle 
with a given orbit decomposition is determined by its structure maps.
We describe a procedure for finding all non-connected quandle structures 
on a disjoint union of subquandles.
\end{abstract}

\textsc{Keywords:} finite quandles, subquandles

\textsc{2000 MSC:} 57M27

\section{\large{\textbf{Introduction}}}

A \textit{quandle} is a set $Q$ with a binary operation $\tr$ satisfying
\newcounter{q}
\begin{list}{(\roman{q})}{\usecounter{q}}
\item{$x\tr x=x$  for every $x\in Q$,}
\item{For every pair $x,y\in Q$, there is a unique $z\in Q$ such that
$x=z\tr y$, and}
\item{For every $x,y,z\in Q$, we have $(x\tr y)\tr z=(x\tr z)\tr(y\tr z)$.}
\end{list}

If $(Q,\tr)$ satisfies (ii) and (iii), $Q$ is a \textit{rack}. Quandles form
a category with morphisms $\phi:Q\to Q'$ defined as maps which preserve the 
quandle operation, i.e.
\[\phi(q_1\tr q_2) = \phi(q_1)\tr'\phi(q_2)\]
where $\tr$ is the quandle operation in $Q$ and $\tr'$ is the quandle
operation in $Q'$. A bijective quandle homomorphism is a quandle isomorphism,
as expected. Axiom (ii) implies that the map $f_b:Q\to Q$ defined by
$f_b(a)=a\tr b$ is bijective for all $b\in Q$; the inverse then defines a 
second operation $a\tr^{-1} b = f_b^{-1}(a)$, called the \textit{dual} 
operation of $Q$; the quandle $(Q,\tr^{-1})$ is the \textit{dual} of $Q$. 
Quandles have been studied in many recent papers; see \cite{FR} for more on 
the history of racks and quandles.

Standard examples of quandles include Alexander quandles, i.e., modules 
$M$ over the ring $\Lambda=\mathbb{Z}[t^{\pm 1}]$ of Laurent polynomials in 
one variable with integer coefficients with quandle operation given by 
\[x\tr y = tx + (1-t) y,\] and groups, which are quandles with quandle
operation given by $x\tr y = y^{-1}xy$. If the group is abelian, the
quandle operation reduces to $x\tr y =x$, and the quandle is 
\textit{trivial}.

Quandles are of interest to topologists since the knot quandle (see \cite{J})
is a complete invariant of knots up to homeomorphism of topological pairs.
Finite quandles are of particular interest as a source of computable knot
invariants such as the counting invariants $|\mathrm{Hom}(K,C)|$ where
$K$ is a knot quandle and $C$ is a finite coloring quandle as well as related
invariants which make use of various quandle cohomology theories (see 
\cite{CGS}, \cite{CJKLS}, etc.)

In this paper, we study the structure of finite quandles in terms of 
subquandles. Our initial goal was to try to find something like a Sylow 
theorem for finite quandles. In section \ref{pd}, we study the 
structure of finite quandles in terms of orbit subquandles. We show 
how to determine all quandles with a given two-subquandle orbit 
decomposition and we discuss how to find quandle structures on a union of 
three or more orbit subquandles.

In section \ref{qm}, we use quandle matrices (see \cite{HN}) to 
study the structure of finite quandles. We describe algorithms for finding the 
orbit decomposition of a finite quandle and for finding quandle structures on 
a disjoint union of $n$ subquandles. This is intended to lay the groundwork 
for the related problem of counting the number of ways of filling in zeroes in 
a quandle presentation matrix to obtain
a finite quandle; it is hoped that a solution to this problem might give new 
insights into the quandle-counting invariants of knots and links studied
in various recent papers (\cite{CJKLS}, \cite{DL}).
\textit{Maple} code for finding rack actions and orbit decompositions
is available at \texttt{www.esotericka.org/quandles}.

\section{\large\textbf{Orbit Decomposition}} \label{pd}

Let $Q$ be a quandle. A \textit{subquandle} $X\subset Q$ is a subset of $Q$ 
which is itself a quandle under $\tr$. Unlike the case of groups, in which the 
intersection of any collection of subgroups is always non-empty (containing 
at least the identity element), a collection of subquandles of a given 
quandle may be pairwise disjoint. Indeed, unlike groups, every subset 
$X\subset Q$ which is closed under $\tr$ is a subquandle: 
if $X$ is closed under $\tr$, the restriction
$f_b|_{X}:X\to X$ is bijective for every $b\in X$, and axiom 
(ii) is satisfied. Since axioms (i) and (iii) are automatic for any subset of 
$Q$, this makes $X$ a subquandle. Thus we have

\begin{lemma} \label{sub}
Let $(Q,\tr)$ be a quandle and $X\subset Q$ a subset. Then $X$ is a subquandle 
iff $X$ is closed under $\tr.$
\end{lemma}

We would like to understand the structure of a quandle in terms of
its subquandles. A quandle which can be written as a union of two disjoint 
subquandles has been called \textit{decomposable} in the literature 
(see \cite{EG}, \cite{G}, \cite{LR}, etc.), and a quandle which is not a 
disjoint union of two subquandles is \textit{indecomposable}. 
The existence of indecomposable quandles follows from the fact that the 
complement of a subquandle is not necessarily a subquandle. However, as 
observed in \cite{LR}, indecomposability of a quandle $Q$ does not imply that 
$Q$ has no subquandles, nor even that the quandle cannot be decomposed as a 
disjoint union of three or more subquandles. 

Indeed, every singleton subset of a quandle is itself a subquandle, though
the analogous statement is not true for non-quandle racks. Thus every quandle 
decomposes in an unhelpful way as a disjoint union of singleton subquandles.
In \cite{LR}, we find the dihedral quandle $R_9$, an
indecomposable quandle which can be written as a disjoint union of three 
isomorphic subquandles. Though this quandle is ``indecomposable,'' since the 
complements of each of the three subquandles are not closed under 
$\tr$, it nevertheless has an internal structure determined by its component 
subquandles -- namely, a Cartesian product of a quandle of order three with 
itself. This is an example of a congruence structure (see \cite{RH}).

\begin{definition} \textup{
Let $Q$ be a quandle and $X\subset Q$ a subquandle. We say that $X$ is 
\textit{complemented in $Q$} or \textit{$Q$-complemented} if $Q\setminus X$ 
is a subquandle of $Q$. Note that since the empty set $\emptyset$ is 
a quandle, every quandle $Q$ is complemented in itself. A quandle $Q$ is 
\textit{complementary} if it has a nonempty $Q$-complemented subquandle. }
\end{definition}

A complementary quandle $Q$ may have subquandles which are not 
$Q$-complemented; indeed, every singleton subset of $Q$ is a subquandle, 
while in general $Q\setminus \{x\}$ is not a subquandle. It is clear from 
the definition of decomposability that a quandle is decomposable iff it is 
complementary.

\cite{J} includes the following definition:

\begin{definition}\textup{
A quandle is \textit{algebraically connected} or just \textit{connected} if 
for every $a,b \in Q$, we have
\[ (\dots ((a \diamond_1 x_1) \diamond_2 x_2) \dots )\diamond_n x_n = b
\]
for some $x_1,\dots, x_n\in Q,\ \diamond_1\dots\diamond_n\in\{\tr,\tr^{-1}\}$. 
The set of all such $b\in Q$ is the \textit{orbit} of $a$.}
\end{definition}

It is well-known that algebraic connectedness coincides with 
indecomposability in the sense defined above, and hence coincides with 
non-complementarity. 

\begin{lemma} \label{intersect}
Let $X$ and $Y$ be subquandles of a quandle $Q$. Then $X\cap Y$ is a 
subquandle. If $X$ and $Y$ are $Q$-complemented, so is $X\cap Y$.
\end{lemma}

\begin{proof}
Let $x,y\in X\cap Y$. Then $x\tr y\in X$ since $X$ is a subquandle,
and $x\tr y\in Y$ since $Y$ is a subquandle. Hence $x\tr y\in X\cap Y$,
and $X\cap Y$ is closed under $\tr$, and $X\cap Y$ is a subquandle by 
lemma \ref{sub}.

Now, suppose $X$ and $Y$ are $Q$-complemented; we must show that 
$Z=Q\setminus (X\cap Y) =(Q\setminus X) \cup (Q\setminus Y)$ is a subquandle. 
Let $x,y\in Z$. If $x$ and $y$ are both in $Q\setminus X$ or both in 
$Q\setminus Y$ then  $x\tr y\in Z$ since $Q\setminus X$ and $Q\setminus Y$
are closed under $\tr$. If $x\in Q\setminus X$ and $y\not\in Q\setminus X$, 
then $y\in X$, which implies $w=x\tr y\in Q\setminus X\subset Z$, since 
otherwise the closure of $X$ under $\tr^{-1}$ would imply $w\tr^{-1} y=x\in X$,
contradicting our choice of $x$. Similarly, $y\in Q\setminus Y$ and 
$x\in Z\setminus(Q\setminus X)$ implies $y\tr x\in Q\setminus Y\subset Z$, 
and $Z$ is closed under $\tr$ as required.
\end{proof}

\begin{theorem} \label{primdec}
Let $Q$ be a finite quandle. Then $Q$ may be written as 
\[Q=Q_1\amalg Q_2 \amalg \dots \amalg Q_n\] 
where every  $Q_i$ is $Q$-complemented and no proper subquandle of any 
$Q_i$ is $Q$-complemented. This decomposition is well-defined up to 
isomorphism; if $Q\cong Q'$, then in the decompositions 
\[ Q=Q_1\amalg \dots \amalg Q_n \quad \mathrm{and} \quad 
Q'=Q_1'\amalg \dots \amalg Q_m',\] we have
then $n=m$ and (after reordering if necessary), $Q_i\cong Q_i'.$
\end{theorem}

\begin{remark}\textup{
The decomposition of a finite quandle into orbits coincides with our notion 
of decomposition into $Q$-complemented subquandles; this follows from the 
observation that the orbits in $Q$ are $Q$-complemented subquandles. 
$Q$-complemented subquandle decomposition then gives us a new perspective on 
the division of $Q$ into disjoint orbits. Indeed, we will see how to 
construct a quandle with a specified list of orbits, when such exists.
See also \cite{EG} proposition 1.17.}
\end{remark}

\begin{proof}
For every $a\in Q$, define $S(a)$ to be the intersection of all 
$Q$-complemented subquandles of $Q$ containing $a$. The collection
$\{S(a) \ | \ a\in Q\}$ is the orbit decomposition of $Q$: each
$S(a)$ is $Q$-complemented, no proper subquandle of any $S(a)$ is 
$Q$-complemented, and $S(a)\cap S(b)\ne \emptyset$ implies $S(a)=S(b)$.
Since the empty quandle is $Q$-complemented, if $Q$ has no nonempty
$Q$-complemented proper subquandles, then $\{S(a) \ | \ a\in Q\}=\{Q\}$; in
any case, $\cup_{a\in Q} S(a)=Q$.

If $\phi:Q\to Q'$ is an isomorphism, then for any subquandle $S\subset Q$
the restriction $\phi|_S$ is an isomorphism onto a subquandle of $Q'$. In 
particular, if $S$ is $Q$-complemented, then $\phi|_{Q\setminus S}$ is also
an isomorphism onto the subquandle $Q'\setminus \phi(S')$. Hence 
$Q'$ has an isomorphic list of $Q'$-complemented subquandles before taking 
intersections, and thus has an isomorphic orbit decomposition.
\end{proof}

\begin{example}\textup{
Let $Q$ be the trivial quandle $T_n=\{1,2,\dots,n\}$ with quandle operation
$i\tr j = i$ for all $i,j\in Q$. Then every singleton subquandle 
$\{i\}\subset Q$ is $Q$-complemented, so the orbit decomposition of $T_n$
is the maximal partition $T_n=\{1\} \amalg \{2\} \amalg \dots \amalg \{n\}$.
}
\end{example}

Before we come to the next theorem, we need a definition.

\begin{definition}\textup{
Let $R$ be a rack and $S$ a quandle. A \textit{rack action} of $R$ on $S$ is
a map from $R$ to the set of automorphisms of $S$, 
$\Phi:R\to \mathrm{Aut}(S)=\{\phi_r:S\to S, \quad r\in R\}$, such that
\[\phi_r (\phi_{r'} (s)) = \phi_{r' \tr r}(\phi_{r}(s))\]
for all $r,r'\in R$ and for all $s\in S$.
}
\end{definition}

\begin{example}\textup{
Let $Q$ be a quandle. Then the set 
$F:Q\to \mathrm{Aut}(Q)=\{f_y:Q\to Q \ | \ f_y(x)=x\tr y\}$
is a rack action of $Q$ on itself, since
\[f_z(f_y(x)) =(x\tr y) \tr z = (x\tr z)\tr (y\tr z) 
= f_{y\tr z}(f_z(x))\]
for all $x,y,z\in Q$.}
\end{example}

\begin{remark}\textup{
Lemma \ref{pdl} is a special case of lemma 1.18 in \cite{EG}.}
\end{remark}

\begin{lemma} \label{pdl}
Let $Q$ and $Q'$ be finite quandles. Then there is a quandle $U=Q\amalg Q'$
iff there are rack actions $F:Q\to \mathrm{Aut}(Q')$ and $G:Q'\to 
\mathrm{Aut}(Q)$ such that the compatibility conditions
\[g_x(a)\tr b = g_{f_b(x)}(a\tr b) \quad \mathrm{and} \quad
f_a(x)\tr y = f_{g_y(a)}(x\tr y) \]
are satisfied for all $a,b\in Q, \ x,y\in Q'$.
\end{lemma}

\begin{proof}
($\Rightarrow$) Let $U=Q\cup Q'$ and let $F,G$ be rack actions. Define 
\[x\tr y =\left\{
\begin{array}{ll}
f_y(x) & x\in Q', y\in Q \\
g_y(x) & x\in Q, y\in Q'.
\end{array}
\right.\] 
Then we assert that $(U,\tr)$ is a quandle. The first quandle axiom is 
satisfied
because $Q$ and $Q'$ are quandles themselves. The second axiom follows 
from the definition of rack action: each $f_a$ and $g_x$ defines a bijection 
of $Q'$ and $Q$ respectively, while the fact that $Q$ and $Q'$ are quandles  
says that each element acts as a bijection on its own subquandle. Hence
the action of $b$ on $U=Q\amalg Q'$ is bijective for each $b\in U$, and axiom 
(ii) is satisfied.

To see that $U$ satisfies (iii), we simply check all the possibilities. 
If $q_1,q_2$ and $q_3$ are in the same subquandle, then (iii) is satisfied.
If $q_1\in Q$ and $q_2,q_3\in Q'$, then 
\[
(q_1\tr q_2)\tr q_3 = g_{q_3}(g_{q_2}(q_1)) = g_{q_2\tr q_3}(g_{q_3}(q_1))
=(q_1\tr q_3)\tr(q_2\tr q_3)
\]
since $G$ is a rack action. Similarly, the fact that $F$ is a rack action
implies that $(q_1\tr q_2)\tr q_3 = (q_1\tr q_3)\tr(q_2\tr q_3)$ when 
$q_1\in Q'$ and $q_2,q_3\in Q$.

If $q_1,q_2\in Q$ and $q_3\in Q'$, then
\[(q_1\tr q_2)\tr q_3 = g_{q_3}(q_1\tr q_2) = g_{q_3}(q_1)\tr g_{q_3}(q_2) 
=(q_1\tr q_3)\tr(q_2\tr q_3)
\]
since $g_{q_3}$ is quandle homomorphism for each $q_3\in Q'$. Similarly
$(q_1\tr q_2)\tr q_3 = (q_1\tr q_3)\tr(q_2\tr q_3)$ when 
$q_1,q_2\in Q'$ and $q_3\in Q$ since each $f_{q_3}$ is a quandle homomorphism.

Finally, if  $q_1,q_3\in Q'$ and $q_2\in Q$, then
\[
(q_1\tr q_2)\tr q_3 = g_{q_3}(f_{q_2}(q_1)) = f_{g_{q_3}(q_2)}(q_1\tr q_3)
=(q_1\tr q_3)\tr(q_2\tr q_3).
\]
Similarly, the compatibility condition ensures that 
$(q_1\tr q_2)\tr q_3 = (q_1\tr q_3)\tr(q_2\tr q_3)$ when 
$q_1,q_3\in Q$ and $q_2\in Q'$.

($\Leftarrow$) If $U=Q\amalg Q'$, then one easily checks that $F:Q\to 
\mathrm{Aut}(Q')$ and $G:Q'\to \mathrm{Aut}(Q)$ are rack actions satisfying 
the compatibility condition. 
\end{proof}

Indeed, as noted in \cite{EG}, lemma \ref{pdl} can be easily generalized to 
obtain:

\begin{theorem} \label{pdra}
Let $Q_1,\dots,Q_n$ be finite quandles. Then there is a quandle
$Q=\amalg_{i=1}^n Q_i$ if there are rack actions 
$\Phi^{i,j}:Q_i\to \mathrm{Aut}(Q_j)$ 
satisfying the compatibility conditions
\[\phi^{k,i}_c(\phi^{j,i}_b(a)) = 
\phi^{k,i}_{\phi^{k,j}_c(b)}(\phi^{k,i}_c(a)) \]
for all $a\in Q_i,$ \ $b\in Q_j$ \ and $c\in Q_k$.
Moreover, $\{Q_1,\dots,Q_n\}$ is the orbit decomposition of $Q$
unless all the rack actions $\Phi^{i,j}$ preserve a $Q_j$-complemented 
subquandle $A\subset Q_j$ for some $Q_j$.
\end{theorem}

\begin{proof}
As above, for each $x\in Q_i, \ y\in Q_j$ define $x\tr y = \phi^{j,i}_y(x)$. 
Then quandle axiom (i) is satisfied automatically since each $Q_i$
is itself a quandle. Axiom (ii) is satisfied since each element
acts on each of the disjoint subquandles $Q_i$ by an automorphism, so the 
overall action is a bijection for each element.

Axiom (iii) is satisfied by the compatibility conditions when the three 
elements are in distinct subquandles or when $a$ and $c$ are in one 
subquandle and $b$ is in another. For example, if $a\in Q_i$, $b\in Q_j$
and $c\in Q_k$ we have
\[
(a\tr b) \tr c = \phi^{j,i}_b(a)\tr c = 
\phi^{k,i}_c(\phi^{j,i}_b(a)) = 
\phi^{k,i}_{\phi^{k,j}_c(b)}(\phi^{k,i}_c(a))
=\phi^{k,i}_{b\tr c}(a\tr c)
=(a\tr c) \tr (b\tr c).
\]

As before, the rack action and 
automorphism requirements satisfy axiom (iii) in the other cases.

Finally, note that removing any $Q_i$ from the list along with the 
corresponding $\Phi^{i,j}$ and $\Phi^{j,i}$ rack actions still defines a
quandle, so the $Q_i$s are each $Q$-complemented. If no subquandle of
any $Q_i$ is preserved by all the actions $\Phi^{i,j}$, then no subquandle
of $Q_i$ is $Q$-complemented and $\{Q_1,\dots,Q_n\}$ is the orbit
decomposition of $Q$.
\end{proof}

\begin{definition}\textup{
Call the maps $\Phi^{i,j}$ the \textit{structure maps} of the quandle $U$
with respect to the decomposition $U=Q_1\amalg Q_2\amalg \dots \amalg Q_n$.}
\end{definition}

\begin{corollary} 
There is a quandle $Q$ with orbits $Q_1,\dots,Q_n$ iff there are compatible
rack actions $\Phi^{i,j}:Q_i\to\mathrm{Aut}(Q_j)$ such that for every
$Q_j$-complemented proper subquandle $S\subset Q_j$ at least one of the
automorphisms $\Phi^{i,j}_a$ for some $a\in Q_i$ does not satisfy 
$\Phi^{i,j}_a(S)=S$.
\end{corollary}

\begin{proof}
If $S\subset Q_j$ is a proper $Q_j$-complemented subquandle such that
every $\Phi^{i,j}_a(S)=S$, then $S$ is an orbit of $Q$.
\end{proof}

\begin{example} \label{trivext}\textup{
Let $Q_1,\dots,Q_n$ be any finite collection of finite quandles, and
define $\Phi^{i,j}:Q_i\to \mathrm{Aut}(Q_j)$ by 
$\phi^{i,j}_x=\mathrm{Id}_{Q_j}$ for all $x\in Q_i$. Then
\[ \phi^{i,j}_x (\phi^{i,j}_y (q)) =q =\phi^{i,j}_{y\tr x} (\phi^{i,j}_x (q))\]
for all $q\in Q_j$, so each $\Phi^{i,j}$ is a rack action. Moreover,
\[\phi^{k,i}_c(\phi^{j,i}_b(a)) = a=
\phi^{k,i}_{\phi^{k,j}_c(b)}(\phi^{k,i}_c(a)) \]
for all $i,j,k$, so $U=Q_1\amalg \dots \amalg Q_n$ is a quandle.
This example shows that there is always at least one quandle structure
on the union of any finite collection of finite quandles. Indeed, if the 
subquandles $Q_1,\dots, Q_n$ are non-complementary, then the orbits
of $U$ are precisely $Q_1,\dots, Q_n$. This example is sometimes called 
the \textit{disjoint union} of the subquandles $Q_1,\dots, Q_n$ (see 
\cite{RH}).}
\end{example}

\begin{example} \label{triv}\textup{
Let $Q_i=\{x_i\}, i =1 \dots n$ be a collection of $n$ singleton quandles.
Then the only possible rack actions by automorphisms of singleton quandles
on other singleton quandles are the identity actions, so the only quandle
structure with orbit decomposition consisting of all singletons must have
$x_i\tr x_j=x_i$ for all $x_i,x_j\in Q$, that is, the trivial quandle $T_n$.}
\end{example}

The observation that a quandle $Q_j\subset Q$ may have subquandles which are
$Q_j$-complemented but not $Q$-complemented implies that quandle may have 
multiple layers of orbit decompositions. Specifically, if 
$U=Q_1\amalg \dots \amalg Q_n$ is the orbit decomposition of $U$, then
each $Q_i$ will have its own orbit decomposition, consisting of multiple
subquandles if $Q_i$ is not connected. Define the \textit{subquandle depth}
of $U$ to be the maximum number $n$ of such layers of decomposition needed
before all remaining orbit decompositions consist of connected subquandles.
This subquandle depth is an invariant of quandle isomorphism type. A connected
quandle has subquandle depth 0; indeed, we may take this as an alternate
definition for ``connected.''

\begin{example}\textup{
The quandle $U$ with quandle matrix (see the next section)
\[
M_U=\left[\begin{array}{rrrr}
1 & 1 & 2 & 2 \\
2 & 2 & 1 & 1 \\
4 & 4 & 3 & 3 \\
3 & 3 & 4 & 4
\end{array}
\right]
\]
has orbit decomposition $U=\{1,2\}\amalg \{3,4\}$. The two orbit subquandles
are trivial and thus have further orbit decompositions $\{1\}\amalg \{2\}$ and
$\{3\}\amalg\{4\}$, so this quandle has subquandle depth 2.}
\end{example}

The structure maps of theorem 5 define a quandle structure on 
$U=Q_1\amalg \dots \amalg Q_n$ with subquandle depth 1. To find all 
non-connected quandle structures on $U$ with subquandles $Q_1,\dots, Q_n$ 
we must consider not only quandle structures with subquandle depth 1 but all 
other possible subquandle depths. To obtain the list of all 
subquandle depth 2 quandle structures on $U$ we must consider all partitions
the set $\{Q_1, \dots, Q_n\}$ into disjoint subsets. Then for each
partition we apply theorem \ref{pdra} to each subset in the partition, 
obtaining the list of all subquandle depth 1 quandles on each subset of the
partition. Applying the theorem again to the resulting new lists of quandles 
yields a set including all subquandle depth 2 quandle structures on $U$. 
Applying this procedure recursively -- that is, for each set $S$ in a partition
of $\{Q_1, \dots, Q_n\}$, consider all the partitions of $S$, etc. --
yields all non-connected quandle structures on $U$ such that each $Q_i$ is 
a subquandle of $U$, since every non-connected quandle structure has some
subquandle depth between 1 and $n$. Note that the sets of quandle structures 
obtained from 
distinct partitions of $\{Q_1, \dots, Q_n\}$ are not disjoint -- the 
structure in which all rack actions are trivial, for example, can be built 
from any partition (or partitioned partition, etc.).

Summarizing, we have

\begin{corollary} \label{pdc}
Every non-connected quandle structure on $U=Q_1\amalg \dots \amalg Q_n$ 
with $Q_i$ subquandles has an orbit decomposition recursively obtainable 
from quandles with orbit decompositions consisting of subsets of 
$Q_1\amalg \dots \amalg Q_n$.
\end{corollary}


Finally, if $\rho:Q\to Q'$ is an isomorphism of quandles where 
$Q=Q_1\amalg Q_2 \amalg \dots \amalg Q_n'$ and
$Q'=Q_1'\amalg Q_2' \amalg \dots \amalg Q_n'$
are the orbit decompositions of $Q$ and $Q'$, then denoting
$\rho|_{Q_i}=\rho_i$, we have
\[ \rho_j (\phi^{i,j}_b(a))= \rho(a\tr b) =\rho(a)\tr \rho(b) =
\phi^{i,j}_{\rho_i(b)} (\rho_j(a)).\] Hence we have

\begin{proposition}
Let $Q_1, \dots, Q_n$ be finite quandles. Then two quandle structures
on the union $U=Q_1\amalg\dots\amalg Q_n$ given by $\Phi^{i,j}$ and
$\Psi^{i,j}$ are isomorphic iff there are automorphisms $\rho_i:Q_i\to Q_i$
such that
\[
\rho_j (\phi^{i,j}_b(a)) = \psi^{i,j}_{\rho_i(b)} (\rho_j(a))
\]
for all $a\in Q_j, b\in Q_i$.
\end{proposition}

\section{\large\textbf{Quandle Matrices and computation}} \label{qm}

The lists of quandles of order $n\le 6$ in \cite{C} and order $n\le 5$ in 
\cite{HN} show that many of the possible quandle structures of small order 
may be
understood as unions of disjoint subquandles. This observation naturally 
raises the question of how many different ways there are for two (or more) 
quandles 
to be put together, that is, how many quandle structures are possible on
$Q\amalg Q'$ which have $Q$ and $Q'$ as subquandles. Theorem \ref{pdra}
gives us an answer, but we need some more convenient notation in order
to permit computations. Quandle matrix notation (see \cite{HN}) provides 
a solution for this problem. 

\begin{definition}\textup{
Let $Q=\{x_1,x_2,\dots, x_n\}$ be a quandle. The \textit{matrix of Q}, 
$M_Q$, is the matrix abstracted from the operation table of $Q$ by
forgetting the $x$s and keeping only the subscripts. That is, we set
$(M_Q)_{ij} =k $ where $x_k=x_i\tr x_j$ in $Q$. }
\end{definition}

Note that quandle axiom (i) enables us to deduce row and column labels 
and hence recover $Q$ from $M_Q$.

\begin{example}\textup{
The Alexander quandle $Q_3=\Lambda_3/(t+1)=\{x_1=0,x_2=1,x_3=2\}$ has 
operation table
\[
\begin{array}{r|rrr}
    & x_1 & x_2 & x_3 \\ \hline
x_1 & x_1 & x_3 & x_2 \\
x_2 & x_3 & x_2 & x_1 \\
x_3 & x_2 & x_1 & x_3
\end{array}
\quad \mathrm{and \ thus} \quad
M_Q=\left[
\begin{array}{rrr}
1 & 3 & 2 \\
3 & 2 & 1 \\
2 & 1 & 3
\end{array}
\right].
\]}
\end{example}

A non-trivial quandle may have an orbit decomposition into trivial
subquandles, and the orbits of a quandle need not be connected.

\begin{example}\textup{
The quandle defined by the quandle matrix
\[\left[\begin{array}{rrrr}
1 & 1 & 1 & 2 \\
2 & 2 & 2 & 3 \\
3 & 3 & 3 & 1 \\
4 & 4 & 4 & 4
\end{array}\right]\]
has orbit decomposition as $T_3\amalg T_1$ with rack actions
$F:T_3\to T_1$ given by $f_1(x)=f_2(x)=f_3(x)=x$ and
$G:T_1\to T_3$ given by $g_4(1)=2,g_4(2)=3,$ and $g_4(3)=1$, that is,
$g_4$ is the permutation $(123)$.}
\end{example}

Just as matrix notation provides a convenient way to specify finite quandles,
we can use column vector notation to represent maps from one finite quandle
to another. Specifically, given a map $\phi:Q\to Q'$ where 
$Q=\{1,2,\dots, n\}$ and $Q'=\{1,2,\dots, m\}$ are quandles given by matrices
$M_Q$ and $M_Q'$, we can represent the map $\phi$ as the $n$-component column 
vector
\[\phi=\left[
\begin{array}{c}
\phi(1) \\
\phi(2) \\
\vdots \\
\phi(n)
\end{array}
\right], \quad \phi(i)\in \{1,\dots,m\}.
\]

Then a rack action $F:Q\to \mathrm{Aut}(Q')$ may be represented as an 
$m\times n$ matrix where the $i$th column is the vector representation of 
$f_i:Q'\to Q'$. Lemma \ref{pdl} then gives us an algorithm for determining 
all quandle structures on $Q\amalg Q'$, namely let $Q=\{1,2,\dots, n\}$ and 
$Q'=\{n+1,\dots,n+m\}$ be quandles with matrices $M_Q$ and $M_{Q'}$ 
respectively. Then
\newcounter{al1}
\begin{list}{(\arabic{al1})}{\usecounter{al1}}
\item{For every $m\times n$ matrix $F$ with columns which are permutations 
of $Q'$, check whether the matrix satisfies the rack action condition
\[F[i,j] = F[M_Q[i,j],i] \quad \forall i,j \in Q\] }
\item{For every $n\times m$ matrix $G$ with columns which are permutations 
of $Q$, check whether the matrix satisfies the rack action condition
\[G[i,j] = G[M_{Q'}[i,j],i] \quad \forall i,j\in Q'\]}
\item{For every pair $F,G$ of such matrices, test the compatibility conditions
\[M[F[i,k],j] = G[M'[i,k],F[G[k,i]]] \quad \forall i,k\in Q, j\in Q' \]
and
\[M'[G[i,k],j] = F[M[i,k],G[F[k,i]]] \quad \forall i,k\in Q', j\in Q. \]}
\item{For every pair $F,G$ which passes steps (1)-(3), the block matrix
\[
\left[\begin{array}{c|c}
M_Q & G \\ \hline
F & M_Q'
\end{array}
\right]
\]
is a quandle matrix.}
\end{list}

Conversely, given a quandle matrix $Q_M$, we can read off the rack actions
by simply interpreting $Q_M$ as a block matrix.

We note that the generalization of this procedure to unions of more than
two quandles does not give all possible quandle structures on the disjoint 
union of three or more subquandles, since this construction yields only 
quandles in which every given $Q_i$ is $Q$-complemented, i.e., quandle 
structures of subquandle depth 1. 

For example, the connected quandle $Q_3\times Q_3$ has matrix
\[
M_{Q_3\times Q_3} =\left[
\begin{array}{ccc|ccc|ccc}
1 & 3 & 2 & 7 & 9 & 8 & 4 & 6 & 5 \\ 
3 & 2 & 1 & 9 & 8 & 7 & 6 & 5 & 4 \\
2 & 1 & 3 & 8 & 7 & 9 & 5 & 4 & 6 \\ \hline
7 & 9 & 8 & 4 & 6 & 5 & 1 & 3 & 2 \\
9 & 8 & 7 & 6 & 5 & 4 & 3 & 2 & 1 \\
8 & 7 & 9 & 5 & 4 & 6 & 2 & 1 & 3 \\ \hline
4 & 6 & 5 & 1 & 3 & 2 & 7 & 9 & 8 \\ 
6 & 5 & 4 & 3 & 2 & 1 & 9 & 8 & 7 \\
5 & 4 & 6 & 2 & 1 & 3 & 8 & 7 & 9
\end{array}
\right],
\]
which has no $Q_3\times Q_3$-complemented subquandles. This quandle is 
isomorphic to the dihedral quandle $R_9$ whose three-subquandle decomposition 
is noted in \cite{LR}. Indeed, we can use the division algorithm to write a 
quandle matrix for $Q\times Q'$ where $|Q|=n$ and $|Q'|=m$ by identifying 
$(x,y)$ with $(x-1)m+y$ for $x=1,\dots, n$ and $y=1,\dots, m$. Then the matrix 
of $Q\times Q'$ is the block matrix
\[\left[
\begin{array}{c|c|c|c}
(q_{11}-1)m + M_{Q'} & (q_{12}-1)m + M_{Q'} & \dots 
& (q_{1m}-1)m + M_{Q'} \\ \hline
(q_{21}-1)m + M_{Q'} & (q_{22}-1)m + M_{Q'} & \dots 
& (q_{2m}-1)m + M_{Q'} \\ \hline
\vdots & \vdots & \ddots & \vdots \\ \hline
(q_{m1}-1)m + M_{Q'} & (q_{m2}-1)m + M_{Q'} & \dots 
& (q_{mm}-1)m + M_{Q'}
\end{array}
\right]\]
where $M_Q=(q_{ij})$ and $M_{Q'}$ are the matrices of $Q$ and $Q'$ 
respectively.

\begin{remark} \label{rem3}\textup{
In the last section, we noted that if a finite quandle $Q$ is a union of 
three or more subquandles, then some quandle structures may have subquandle
depth greater than 1, since $x\in Q_1$ need not imply $x\tr y \in Q_1$. For 
example, the quandle $Q$ below is a union of three subquandles $Q_1=\{1,2\}, 
Q_2=\{3,4\}$ and $Q_3=\{5,6\}$, and indeed 
there is an apparent block-matrix decomposition. However, because 
$Q_1$ is not $Q$-complemented, there is no rack action $\Phi^{2,1}:
Q_2\to\mathrm{Aut}(Q_1)$, for example.
\[
M_Q= \left[
\begin{array}{cc|cc|cc}
1 &   1 &   2 &   2 &   1  &  1 \\
2 &   2 &   5 &   5 &   2  &  2 \\ \hline
3 &   3 &   3 &   3 &   3  &  3 \\
4 &   4 &   4 &   4 &   4  &  4 \\ \hline
5 &   5 &   1 &   1 &   5  &  5 \\
6 &   6 &   6 &   6 &   6  &  6
\end{array}
\right]
\]
To construct this subquandle depth 2 quandle from $Q_1, Q_2$ and $Q_3$ we 
must first put together $Q'=Q_1\amalg Q_3\cong T_4$, then find structure 
maps for $Q=Q'\amalg Q_2$.}
\end{remark}

\textit{Maple} programs for finding rack actions and orbit decompositions
of finite quandles represented by matrices are available in the file 
\texttt{quandles-maple.txt} at \texttt{www.esotericka.org/quandles}.


\begin{thebibliography}{00}

\bibitem{EG}  N. Andruskiewitsch and M. Gra\~{n}a. 
From racks to pointed Hopf algebras.  \textit{Adv. Math.}  \textbf{178} 
(2003) 177-243.

\bibitem{C} J. S. Carter, S. Kamada and M. Saito. 
\textit{Surfaces in 4-space.}
Encyclopaedia of Mathematical Sciences, 142. Low-Dimensional Topology, III.
Springer-Verlag, Berlin,  2004.

\bibitem{CJKLS} J. S. Carter, D. Jelsovsky, S. Kamada, L. Langford and M. 
Saito. 
Quandle cohomology and state-sum invariants of knotted curves and surfaces.  
\textit{Trans. Amer. Math. Soc.}  \textbf{355} (2003) 3947-3989.

\bibitem{CGS} J. S. Carter, M. Elhamdadi, M. Gra\~{n}a and M. Saito. 
Cocycle Knot Invariants from Quandle Modules and Generalized Quandle 
Cohomology.  \textit{Osaka J. Math.}  \textbf{42}  (2005) 499-541.

\bibitem{DL} F. M. Dion\'{i}sio and P. Lopes. 
Quandles at finite temperatures. II.  
\textit{J. Knot Theory Ramifications} \textbf{12} (2003) 1041-1092.

\bibitem{FR} R. Fenn and C. Rourke. 
 Racks and links in codimension two. 
\textit{J. Knot Theory Ramifications}  \textbf{1}  (1992) 343-406.

\bibitem{G} P. Etingof and M. Gra\~{n}a. 
 On rack cohomology. \textit{J. Pure Appl. Algebra} \textbf{177} (2003) 49-59.

\bibitem{MN} R. Henderson, T. Macedo and S. Nelson
Symbolic Computation with finite quandles. \textit{J. Symbolic Comput.} 
\textbf{41} (2006) 811-817.

\bibitem{HN} B. Ho and S. Nelson. 
Matrices and Finite Quandles.
\textit{Homology Homotopy Appl.} \textbf{7} (2005) 197-208.

\bibitem{J} D. Joyce. A classifying invariant of knots, the knot quandle.  
\textit{J. Pure Appl. Algebra}  \textbf{23}  (1982)  37-65.

\bibitem{LR} P. Lopes and D. Roseman. 
On finite racks and quandles.  
\textit{Comm. Algebra} \textbf{34} (2006) 371-406.

\bibitem{RH} H. Ryder. 
\textit{The Structure of Racks}. Ph.D. Dissertation, U. Warwick. Available 
online at \texttt{www.esotericka.org/quandles}

\end{thebibliography}
\end{document}